\documentclass[runningheads]{llncs}

\usepackage{cite}
\usepackage{amssymb}
\usepackage{rotating}
\usepackage{amsmath}
\usepackage{graphicx}
\usepackage{color}
\usepackage[utf8]{inputenc}
\usepackage[american]{babel}
\usepackage{enumitem}
\usepackage{caption}
\usepackage{nicefrac}
\usepackage{booktabs}
\usepackage{hyperref}
\usepackage{todonotes}
\usepackage{subcaption}

\graphicspath{{./figures/}}

\begin{document}
\setlength{\tabcolsep}{4pt}
\title{Idle vehicle repositioning for dynamic ride-sharing}

\author{Martin Pouls\inst{1}\orcidID{0000-0002-9258-719X}\and
Anne Meyer\inst{2}\orcidID{0000-0001-6380-1348}\and
Nitin Ahuja\inst{3}}
\authorrunning{M. Pouls et al.}

\institute{FZI Research Center for Information Technology, 76131 Karlsruhe, Germany \\
\email{pouls@fzi.de}\\ \and
TU Dortmund University, 44221 Dortmund, Germany \and
PTV Group, 76131 Karlsruhe, Germany}

\maketitle  

\begin{abstract}
In dynamic ride-sharing systems, intelligent repositioning of idle vehicles enables service providers to maximize vehicle utilization and minimize request rejection rates as well as customer waiting times. In current practice, this task is often performed decentrally by individual drivers.
We present a centralized approach to idle vehicle repositioning in the form of a forecast-driven repositioning algorithm. The core part of our approach is a novel mixed-integer programming model that aims to maximize coverage of forecasted demand while minimizing travel times for repositioning movements. This model is embedded into a planning service also encompassing other relevant tasks such as vehicle dispatching. 
We evaluate our approach through extensive simulation studies on real-world datasets from Hamburg, New York City, and Manhattan. We test our forecast-driven repositioning approach under a perfect demand forecast as well as a naive forecast and compare it to a reactive strategy. The results show that our algorithm is suitable for real-time usage even in large-scale scenarios. Compared to the reactive algorithm, rejection rates of trip requests are decreased by an average of 2.5 percentage points and customer waiting times see an average reduction of $13.2 \%$.

\end{abstract}
\section{Introduction}
\label{sec:introduction}
While the popularity of mobility-on-demand (MOD) services such as Uber and Lyft has increased significantly in recent years, this growth has also lead to increased traffic congestion \cite{sfcta_2018}. Several cities have identified this issue and some have even taken countermeasures \cite{npr_2018}. One way to tackle this problem is the increased usage of dynamic ride-sharing services such as UberPool or MOIA. In these services, multiple passengers with different destinations share a vehicle. This way, one maintains the flexibility of MOD services compared to traditional public transport, while at the same time improving vehicle utilization.

\begin{figure}[htbp]
\centering
\includegraphics[width=0.8\textwidth]{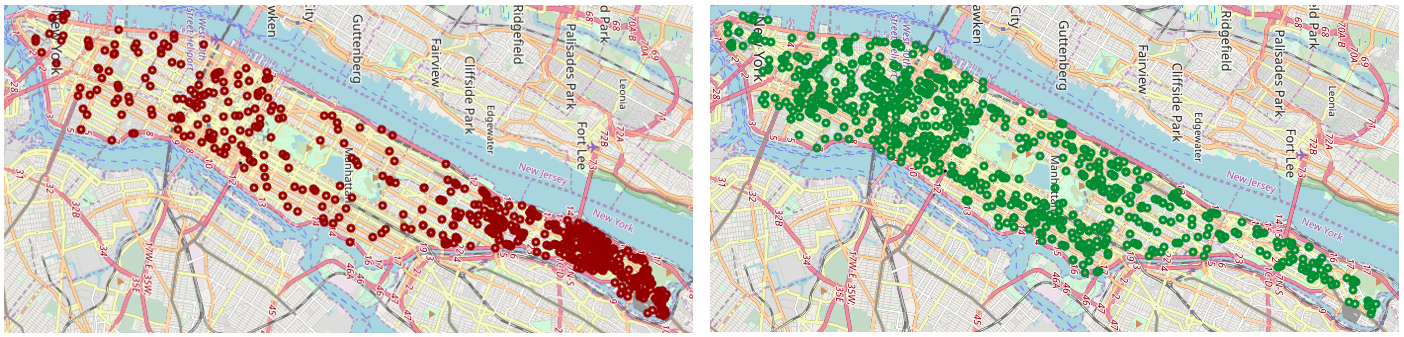}\caption{Vehicle positions without (left) and with (right) repositioning.}
\label{fig:motivation}
\end{figure}

Planning problems regarding MOD services in general and dynamic ride-sharing, in particular, have generated significant research attention. Most works focus on the vehicle routing aspect, i.e. solving the dynamic dial-a-ride-problem arising in these applications \cite{alonso-mora_-demand_2017, ma_dynamic_2018}. In this work, we focus on the idle vehicle repositioning problem, i.e. the problem of sending idle vehicles to a suitable location in anticipation of future demand. In many practical applications with self-employed drivers, this problem is currently solved decentrally by incentivizing drivers to reposition towards areas with low vehicle supply. For instance, Uber employs so-called "surge pricing" which raises prices in areas with excess demand and thereby offers increased revenue opportunities to drivers \cite{uber_2020}. We propose the usage of a central repositioning strategy that may improve system performance in use cases with a central fleet operator.
In general, the overall performance of a ride-sharing system may be impacted significantly by suitable repositioning algorithms. Figure \ref{fig:motivation} illustrates this fact by comparing vehicle positions in scenarios without and with repositioning.
Without repositioning vehicles become stuck in low-demand areas. In turn, requests in other areas are rejected due to a lack of nearby vehicles. This is due to the assumption of a maximum waiting time for customers in dial-a-ride problems. A vehicle must reach the customer within this time frame, otherwise the customer is rejected. Thus, vehicles in low-demand areas cannot reach many new trip requests in time and are consequently rarely assigned a new tour. This phenomenon may be avoided by using a suitable repositioning mechanism.

The main contribution of this paper is the introduction of a novel forecast-driven repositioning algorithm. In particular, we propose a mixed-integer programming (MIP) model that maximizes the coverage of forecasted demand while minimizing driving times for repositioning movement. In contrast to prior works, our model utilizes a realistic demand forecast conforming to state-of-the-art forecasting models \cite{liao_large-scale_2018, yao_deep_2018} and does not assume any further information regarding probability distributions of trip requests. The solution approach is embedded into a planning service for dynamic ride-sharing applications and evaluated on real-world taxi datasets from Hamburg, New York City and Manhattan through extensive simulation studies. The results show that our model can be used in real-time even on large-scale instances with up to 20,000 trip requests per hour. Compared to a reactive approach, our algorithm reduces rejection rates by an average of 2.5 percentage points and customer waiting times by $13.2 \%$.

The remainder of this work is organized as follows. Section \ref{sec:related_work} gives an overview of related work regarding idle vehicle repositioning. In Section \ref{sec:system_overview} we briefly describe our planning service and the simulation that is used for evaluations. Our repositioning approach is detailed in Section \ref{sec:repositioning_approaches}. Finally, Section \ref{sec:computational_results} presents our computational results and Section \ref{sec:outlook} summarizes our findings and gives some possible directions for future work.
\section{Related work}
\label{sec:related_work}
There has recently been an influx of papers dealing with repositioning in the context of MOD services. For bike- and car-sharing applications, several repositioning strategies have been proposed (e.g. \cite{chemla_bike_2013, weikl_relocation_2013}). However, all of these approaches differ from the problem at hand either due to their station-based nature or due to the missing consideration of ride-sharing. In the context of classical taxi services, some works consider repositioning and propose the usage of historical GPS data to identify potentially profitable regions \cite{li_hunting_2011, pfoser_towards_2011}. In contrast to our work, these papers do not consider ride-sharing and have different objective functions as they view the problem decentrally from the perspective of individual drivers.
Idle vehicle repositioning has also been modeled using queuing-based methods \cite{sayarshad_non-myopic_2017, braverman_empty-car_2019}. Although these works show that their approaches yield improvements compared to myopic strategies, they have not been evaluated on realistic large-scale scenarios. Moreover, they tend to be limited in the extent of the covered area and the spatial granularity of decisions.

To the best of our knowledge, there are only few papers considering idle vehicle repositioning in the context of large-scale dynamic ride-sharing applications. 
One approach is a reactive repositioning scheme with the idea of sending idle vehicles to the pickup locations of rejected trip requests \cite{alonso-mora_-demand_2017}. Idle vehicles are matched to rejected requests while minimizing the travel time for repositioning movements. 
In a follow-up paper \cite{alonso-mora_predictive_2017}, the same authors present a more refined sampling-based approach. They include anticipated trip requests in their vehicle routes that are served with a low priority. The authors show that this approach leads to reduced waiting times and in-car travel delays compared to the reactive repositioning from \cite{alonso-mora_-demand_2017}. However, no significant improvement in the number of served trip requests was made.
Another paper presents two simple approaches in which vehicles reposition according to historical pickup probabilities \cite{jung_large-scale_2019}. Vehicles either move to a zone for roaming or to a depot. The probability of selecting a zone or depot is proportional to historical pickup probabilities. The authors compare these approaches to a setting without repositioning and evaluate them with a simulation scenario based on New York taxi data with approximately 145,000 trip requests over 8 hours. Both repositioning algorithms improve the request acceptance rate at the cost of an increase in deadheading time. In contrast to our work, the authors do not consider detailed information about supply and demand. In particular, neither the current configuration of the vehicle fleet nor the total demand is considered during repositioning.

In this work, we view repositioning as an independent problem. As seen from \cite{alonso-mora_predictive_2017}, repositioning can also be treated as an integrated decision during vehicle routing. In that case, the problem may be viewed as a vehicle routing problem with stochastic customers. A variety of solution approaches have been presented for this problem, for a review see \cite{ritzinger_survey_2016}. However, none of these approaches have been studied on large-scale scenarios and they often assume the presence of detailed information regarding trip request distributions that is not available in many practical settings.
\section{System Overview}
\label{sec:system_overview}
\begin{figure}[htbp]
\includegraphics[width =\linewidth]{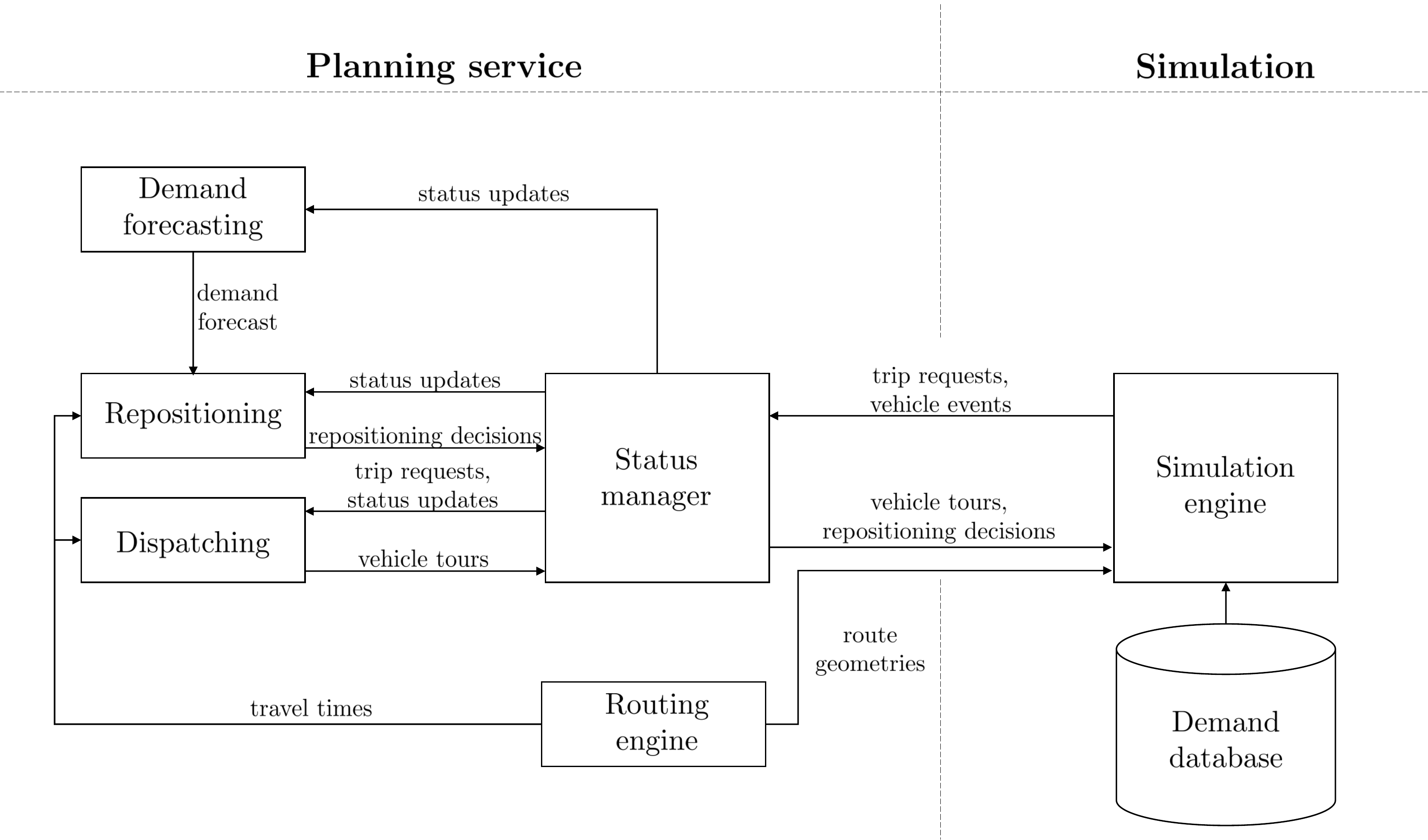}
\caption{Planning service, simulation and relevant communication.}
\label{fig:system_central}
\end{figure}

Our forecast-driven repositioning algorithm is embedded into a planning service for dynamic ride-sharing applications encompassing all components regarding dispatching, repositioning, demand forecasting, and routing. For evaluation, the planning service is coupled with a simulation that emulates relevant real-world events generated by customers and vehicles. The resulting overall system and the communication between components is depicted in Figure \ref{fig:system_central}.

By making a strict separation from the simulation, the planning service could theoretically be directly transferred to a real-world use case.

\subsection{Planning service}
\label{sec:system_overview:dispatching}
All planning and forecasting functionality is consolidated in the planning service, which consists of separate decoupled modules. All external communication takes place via the central status manager that also maintains the current status of all vehicles and requests and provides this information to other modules.

The dispatching module contains functionalities for vehicle routing. Essentially, it solves a dynamic dial-a-ride-problem with the typical constraints: capacity, waiting time, and ride time \cite{cordeau_dial--ride_2007}. We utilize a simple insertion heuristic inspired by \cite{ma_real-time_2015} which tries to dynamically insert each incoming request into the current routing plan. If no feasible insertion is found, the request is rejected. The repositioning module implements our approach for idle vehicle repositioning which is detailed in Section \ref{sec:repositioning_approaches}.

Two further components are needed to provide the necessary input for planning. The routing engine operates on OpenStreetMap (OSM) data and provides travel times to the other modules. The demand forecasting module outputs the forecasted number of trip requests given a discrete set of areas (e.g. a grid-based partitioning of the map) and a forecast horizon (e.g. 30 minutes). This format conforms to state of the art demand forecasting methods \cite{liao_large-scale_2018, yao_deep_2018} that could be integrated in this system. However, the forecasting methodology itself is not the focus of this work.

\subsection{Simulation}
\label{sec:system_overview:sim}
The simulation operates on a database containing trip requests. Each trip request consists of the request time, origin and destination coordinates, and the number of passengers. The request may be obtained from actual taxi services or be derived from other sources such as public transport data or traffic simulations.

These trip requests are replayed by the simulation engine and sent to the planning service. If the trip request is accepted, it is assigned to a vehicle and the simulation receives an updated vehicle tour. It operates on this tour and simulates all relevant events such as arrival and departure at stops. Furthermore, the simulation emulates real-world GPS tracking by regularly sending position updates to the planning backend. The necessary information for these updates is obtained from a routing engine working on the road network.

\section{Repositioning approaches}
\label{sec:repositioning_approaches}
In the following, we present our forecast-driven repositioning algorithm (FDR) as well as a simple reactive strategy (REACT) intended as a benchmark.

\subsection{Forecast-driven repositioning}
\label{sec:forecast}
Our algorithm works with a demand forecast that provides the anticipated number of trip requests for a set of areas and a forecast horizon. The core part of FDR is a MIP model (FDR-M) in which we aim to maximize the coverage of forecasted demand by intelligently repositioning idle vehicles while minimizing the travel times for these movements. We assume that vehicles may cover trip requests near their current location as they can reach these requests within the maximum waiting time. Our model also takes the current state of vehicles into account and reflects the fact that vehicles may serve multiple requests at once.
Model FDR-M takes decisions on an aggregated level. As an output, it determines the number of vehicles relocated between specific areas. The model is embedded into a planning process that provides the necessary inputs and translates the model output into actual repositioning assignments.
In the following, we will first present this planning process and subsequently introduce the model itself with the necessary notation.

\subsubsection{Planning process}
\label{sec:forecast:planning}
Model FDR-M is embedded into a rolling horizon planning process which is triggered at regular intervals $f$ (e.g. every 3 minutes). The main steps of the planning process are as follows:
\begin{enumerate}
\item Obtain an up-to-date demand forecast.
\item Solve repositioning model FDR-M.
\item Determine actual repositioning assignments.
\end{enumerate}
In the last step, we determine specific vehicles and targets for repositioning based on the aggregated output of FDR-M. Assume the model decides to reposition $x_{ij}$ vehicles from area $i$ to area $j$. Given that $j \in A$ represents an area, we have a set of feasible repositioning targets $P_j$ in $j$. These are specific points to which we may send a vehicle. They are determined from prior trip requests, i.e. each past pickup location is a feasible repositioning target.
We sample $x_{ij}$ of these targets and subsequently greedily send the closest idle vehicle to each target.

\subsubsection{Notation}
\label{sec:forecast:notation}
\begin{figure}[htbp]
    \centering
    \includegraphics[width=0.45\textwidth]{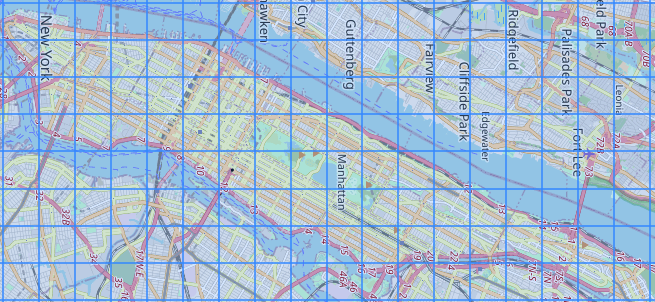}\caption{Grid-based partitioning into areas (1x1 km).}
    \label{fig:grid}
    \end{figure}

Relevant notation for FDR is summarized in Table \ref{tab:notation}.
$K$ denotes the set of all vehicles. This set may be further subdivided into idle vehicles $K^{id}$, vehicles serving a tour $K^{t}$ and vehicles on a repositioning trip $K^{r}$. In addition, we assume a set of areas $A$. For the remainder of this work, we utilize a partitioning of the region under study into grid cells (1x1 km) as seen in Figure \ref{fig:grid}. 
For the purpose of travel time calculations, we assume that these areas are represented by their center and determine a travel time $t_{ij}$ between the centers of $i$ and $j$. We may now further divide the sets of vehicles by area as $K^{id}_a$, $K^t_a$, $K^r_a$. Note that the sets $K^{id}_a$ and $K^t_a$ contain those vehicles currently situated in area $a$. $K^r_a$ on the other hand consists of vehicles currently repositioning towards $a$. Vehicles may only be repositioned to valid target areas $a \in A_r \subseteq A$. In practical applications $A_r$ might be determined based on suitable waiting spots for vehicles. In this study, we limit $A_r$ to areas with at least one prior pickup as we sample specific repositioning targets from past pickup locations. Our model works on a demand forecast denoted as $\hat{d}_{a}, a \in A$ over the forecast horizon $h$ (e.g. 30 minutes). This forecast gives us the predicted number of trip requests originating in $a$ within the horizon $h$. The areas that a vehicle may cover from its current area $a \in A$ are controlled by $t^c$ which corresponds to the maximum waiting time of a customer. We refer to an area $j$ as reachable from $i$ ($r_{ij} = 1$), if $t_{ij} \leq t^c$. The assumption is that in this case a vehicle located at $i$ may reach a request in $j$ on time. Vehicles may serve multiple trip requests over the forecast horizon $h$. This aspect is included in the model by parameter $p$ which corresponds to the assumed number of requests that an initially idle vehicle can serve. Vehicles currently serving a tour may also cover future demand. However, given the fact that they are already partially occupied with their current tour, their provided coverage is discounted by a factor $\alpha < 1$. Both $\alpha$ and $p$ should be determined depending on the specific scenario under study. They may for instance be derived from historical data. A possible extension would be to determine these parameters adaptively and vehicle-specific.
\begin{table}[t]
    \caption{Notation for forecast-driven repositioning.}
    \label{tab:notation}
    \centering
    \begin{tabular}{l l}
    \toprule
    \multicolumn{2}{l}{\textbf{Sets}}\\
    \midrule
    $K$ & Vehicles \\
    $K^{id} | K^t | K^r$ & Idle / touring / repositioning vehicles \\
    $A$ & Areas \\
    $A_r$ & Valid target areas for repositioning \\
    $K^{id}_a | K^t_a | K^r_a$ & Idle / touring / repositioning vehicles per area $a \in A$ \\
    $P_a$ & Set of feasible repositioning targets in $a \in A$ \\ 
    \midrule
    \multicolumn{2}{l}{\textbf{Parameters}}\\
    \midrule
    $t_{ij}$ & Travel time from $i \in A$ to $j \in A$ \\
    $\hat{d}_{a}$ & Demand forecast for $a \in A$ \\
    $h$ & Forecast horizon \\
    $p$ & Served requests over the given forecast horizon $h$ \\
    $t^c$ & Coverage radius \\
    $r_{ij}$ & Reachability indication $\in \{0,1\}$; $r_{ij} = 1$, if $t_{ij} \leq t^c$ \\
    $\alpha$ & Discount factor for touring vehicles $k \in K^t$ \\
    $\beta$ & Factor for travel times incurred by coverage \\
    \bottomrule
    \end{tabular}
    \end{table}

\subsubsection{Mixed-integer formulation}
\label{sec:forecast:mip}
Model FDR-M is given in equations (\ref{eq:mip:obj}) - (\ref{eq:mip:domain_c}). Variables $c_{ij} | i, j \in A$ denote the coverage that is provided by vehicles in $i$ for forecasted demand in $j$. This corresponds to the assumed amount of forecasted trip requests in $j$ that would be served by vehicles from $i$. $c_{ij}$ may take fractional values as forecasted demand $\hat{d}_{j}$ may be fractional and can be covered from multiple origin areas $i \in A$. Variables $x_{ij}| i, j \in A$ denote the number of vehicles repositioned from $i$ to $j$.
The hierarchical objective function~(\ref{eq:mip:obj}) follows three goals which are reflected in the terms of the objective function:
\begin{enumerate}
\item Maximize the sum of covered demand, weighted by forecasted demand.
\item Minimize the number of repositioning movements.
\item Minimize weighted travel times. 
\end{enumerate}
Objective precedence is ensured by weights $W_1 >> W_2 >> 1$. The primary objective is to maximize the acceptance rate of future requests by covering predicted demand. Empirically, it has proven beneficial to prioritize coverage in high demand areas. Therefore, we add weights corresponding to the forecasted demand $\hat{d}_j$ in the covered area $j \in A$. The secondary objective stems from the operational concern that we want to move as few vehicles as possible. Particularly, we do not want to move any vehicles at all, if the current fleet configuration can cover all forecasted demand. The tertiary objective ensures that overall travel times are minimized and leads to suitable vehicles being selected for repositioning. Two travel time factors are considered. On the one hand, travel times are attached to $x_{ij}$ variables as repositioning movements incur the movement of empty vehicles. Additionally, we consider anticipated travel times attached to $c_{ij}$ variables. The assumption is that a vehicle located at $i \in A$ will have to move to $j \in A$ when a request arises. These anticipated travel times are penalized by a factor $\beta \geq 1$ which rewards moving vehicles closer to predicted demand. This tends to be beneficial as it reduces customer waiting times and improves vehicle utilization.

\begin{align}
    & \text{(FDR-M)}& \sum_{i \in A} \sum_{j \in A} & W_1 \cdot \hat{d}_j \cdot c_{ij} - \sum_{i \in A} \sum_{j \in A} W_2 \cdot x_{ij} && \label{eq:mip:obj} \\
    && - \sum_{i \in A} \sum_{j \in A} & (x_{ij} + \beta \cdot c_{ij}) t_{ij} \rightarrow \text{max} && \notag \\
    & \text{s.t.}&\sum_{j \in A} x_{ij} &\leq |K^{id}_i| && i \in A \label{eq:mip:max_idle}\\
    && \sum_{j \in A} c_{ji} & \leq \hat{d}_i && i \in A \label{eq:mip:max_cover}\\
    &&  \sum_{j \in A}c_{ij} & \leq p \cdot (\sum_{j \in A} x_{ji} + |K^r_i| + \alpha \cdot |K^t_i|) && i \in A \label{eq:mip:max_available}\\
    && c_{ij} & = 0 && i, j \in A, r_{ij} = 0 \label{eq:mip:reach}\\
    && x_{ij} & = 0 && i \in A, j \notin A_{r} \label{eq:mip:target}\\
    && x_{ij} & \in \mathbb{N}^+_0 && i, j \in A \label{eq:mip:domain_x}\\
    && c_{ij} & \in \mathbb{R}^+_0 && i, j \in A \label{eq:mip:domain_c}
\end{align} 

Constraints (\ref{eq:mip:max_idle}) guarantee that the number of vehicles repositioned from $i \in A$ does not exceed the number of idle vehicles. Constraints (\ref{eq:mip:max_cover}) ensure that the maximum provided coverage for a given area is capped by the forecasted demand. Inversely, Constraints (\ref{eq:mip:max_available}) limit the provided coverage from area $i$ to the maximum available coverage. This maximum available coverage for area $i$ is calculated based on the available vehicles and includes the assigned vehicles $x_{ji}$. In addition, provided coverage is based on vehicles repositioning to $i$ or currently on a tour in $i$. The latter ones are discounted by factor $\alpha \leq 1$. The available coverage in trip requests is obtained by multiplying the available vehicles with the assumed number of served requests per period $p$. Reachability of covered areas is ensured by Constraints (\ref{eq:mip:reach}) while Constraints (\ref{eq:mip:target}) limit repositioning movements to valid targets. Variable domains are given by Constraints (\ref{eq:mip:domain_x}) and (\ref{eq:mip:domain_c}).

\subsection{Reactive repositioning}
\label{sec:reactive}
As a benchmark, we implement a reactive approach (REACT). The algorithm is an adapted version of the reactive repositioning algorithm presented by the authors of \cite{alonso-mora_-demand_2017}. We modify their approach to reflect the fact that we process trip requests individually whereas they work with batches. Therefore, after rejecting a request, we may also directly reposition an idle vehicle. Given a rejected request $r$ an its pickup location $p_r$, we greedily reposition the nearest idle vehicle to $p_r$, i.e. the vehicle with the shortest travel time from its current position to $p_r$.

\section{Computational results}
\label{sec:computational_results}

\subsection{Data and setup}
\label{sec:computational_results:data}
We evaluate our repositioning algorithms on three real-world taxi datasets from Hamburg\footnote{Provided by PTV Group, Haid-und-Neu-Str. 15, 76131 Karlsruhe} (HH), New York\footnote{https://www1.nyc.gov/site/tlc/about/tlc-trip-record-data.page} (NYC) and Manhattan (MANH). The latter is created from the NYC dataset by limiting it to trips within Manhattan. These datasets contain the pickup time, pickup location, dropoff location, and number of passengers for historical taxi trips. We filter the original data by eliminating obvious outliers and erroneous records. As a routing engine, we use RoutingKit \cite{dibbelt_customizable_2015} which operates on OpenStreetMap extracts covering the respective areas under study. Gurobi 8.1.0 serves as a MIP solver for model FDR-M. All experiments were run on the same machine with an Intel i7-6600U CPU and 20 GB of RAM.

\subsection{Scenarios and parameter settings}
For each dataset, we run simulations covering two separate temporal scenarios: a Wednesday and a Sunday. 
\begin{table}
    \caption{Temporal scenarios with dates and number of trip requests.}
    \label{tab:time_periods}
    \centering
    \begin{tabular}{l l l l l l l}
    \toprule
    & \multicolumn{2}{l}{\textbf{HH}} & \multicolumn{2}{l}{\textbf{NYC}} & \multicolumn{2}{l}{\textbf{MANH}}\\
    \cmidrule (lr){2-7}
    & \textbf{date} & \textbf{req.} & \textbf{date} & \textbf{req.} & \textbf{date} & \textbf{req.}\\
    \cmidrule(lr){2-3}\cmidrule(lr){4-5}\cmidrule(lr){6-7}
    Wed. & 20.03.2019 & 13,556 & 16.03.2016 & 376,526 & 16.03.2016 & 297,457 \\
    Sun. & 24.03.2019 & 10,669 & 20.03.2016 & 368,508 & 20.03.2016 & 269,346 \\
    \bottomrule
    \end{tabular}
    \end{table}
This was done to evaluate the algorithm performance under different demand pattern as the spatial and temporal distribution of trip requests varies between weekdays and weekends. The precise dates and the number of trip requests are shown in Table \ref{tab:time_periods}.

We also vary the size of the vehicle fleet. We determined a base number of vehicles from preliminary testing (HH -- 90, NYC -- 1300, MANH -- 900). With this base number, the fleet should be able to service around $95\%$ of all trip requests. We then create scenarios with vehicle factors of $0.8$, $0.9$, $1.0$, $1.1$ and $1.2$ where the actual number of vehicles per dataset is obtained by multiplying the base number with the vehicle factor. Combined with the temporal settings, we end up with 10 scenarios per dataset. Each of these scenarios is run with three different repositioning modes: no repositioning (NONE), reactive repositioning as described in Section \ref{sec:reactive} (REACT) and forecast-driven repositioning from Section \ref{sec:forecast:mip} (FDR). The latter is run with two different demand forecasts: 1. a perfect demand forecast (FDR (P)) and 2. a naive demand forecast (FDR (N)). This naive demand forecast assumes that demand stays constant, i.e. the forecasted demand over the next horizon $h$ is equal to the observed demand within the previous horizon $h$. 

Parameters concerning FDR were determined based on preliminary results and historical data and are summarized in Table \ref{tab:params}.
\begin{table}
\caption{FDR parameter settings.}
\label{tab:params}
\centering
\begin{tabular}{l l l l l l}
\toprule
\textbf{Description}& \textbf{All} & \textbf{HH} & \textbf{NYC} & \textbf{MANH} \\
\midrule
Forecast horizon ($h$) & 30 min & & &\\
Repositioning frequency ($f$) & 3 min & & & \\
Coverage travel time weight ($\beta$) & 1.05 & & & \\
Active vehicle factor ($\alpha$) & 0.7 & & & \\
Grid cell size & 1x1 km & & & \\
Objective weight for total coverage ($W_1$) & 1000 & & & \\
Objective weight for vehicle movements ($W_2$) & 10 & & & \\
Requests per forecast horizon ($p$) & & 5  & 8 & 9 \\
Coverage limit ($t^c$) & & 8 min & 4 min & 4 min \\
\bottomrule
\end{tabular}
\end{table}
Some parameters need to be determined dataset-specific, due to the significant differences in covered area and in demand density between the three datasets. For instance, in Manhattan a single vehicle may serve more requests in 30 minutes than in Hamburg due to much denser demand and a smaller covered area. The same factors also lead to the longer maximum waiting time for HH compared to MANH and NYC. All simulations are run with a warm-up time of 6 hours, i.e. if we are evaluating the 20.03.2019, the simulation actually starts at 19.03.2019 18:00 and the first 6 hours are not included in the gathered statistics.

\subsection{Algorithm performance with a perfect demand forecast}
In this section, we first evaluate the performance of FDR under a perfect demand forecast. Section \ref{sec:computational_results:naive} compares these results to a setting with a naive forecast.

\subsubsection{Running times}
Our forecast-driven repositioning is real-time capable even on large problem instances. The average running times for one iteration of algorithm FDR was 475 \text{ms} (HH), 1938 \text{ms} (NYC) and 138 \text{ms} (MANH). Given that FDR is run once every 3 minutes, this running time is unproblematic. Including all other tasks such as dispatching and simulation, the average total running time for one scenario run was 8.7 \text{min} (HH), 315.9 \text{min} (NYC) and 138.4 \text{min} (MANH), a substantial speed-up over the simulated real-time equivalent of 1440 \text{min}.

\subsubsection{Rejection rates and vehicle travel times}
\begin{figure}[htbp]
\centering
\includegraphics[width=0.90\textwidth]{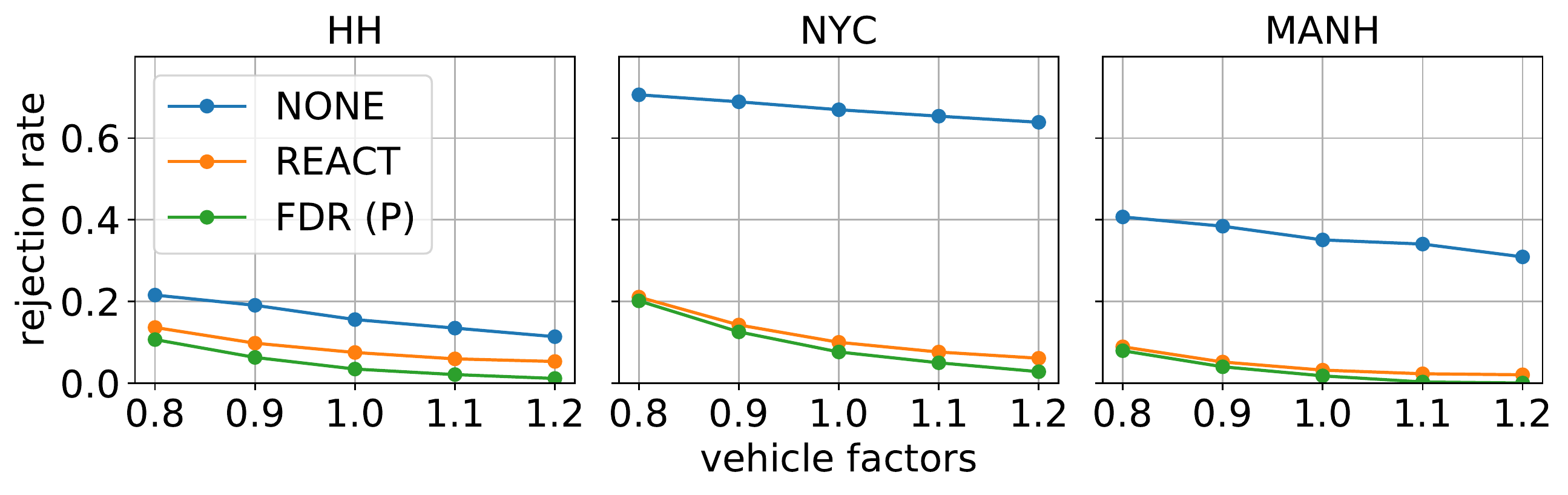}\caption{Average request rejection rates.}
\label{fig:rejection}
\end{figure}

Figure \ref{fig:rejection} shows the request rejection rates, i.e. the fraction of requests that could not be served, depending on the dataset and fleet size factor. Across all scenarios, FDR yields the best results. The average improvement compared to REACT was 3.7, 2.2 and 1.5 percentage points for HH, NYC, and MANH respectively.
The improvement varies substantially between the datasets. We believe this is mainly due to the geographical distribution of requests and overall vehicle utilization. For instance, REACT works remarkably well on the Manhattan scenario where most trip requests occur in downtown Manhattan. In case of the other two datasets, the improvement in rejection rates is more significant.
One trend across all datasets is that the difference between FDR and REACT grows as the number of vehicles is increased. FDR is better suited to exploit larger fleet sizes where almost all trip requests may be served, even ones in remote areas. In case of small vehicle fleets, the complete fleet may be occupied during peak hours, therefore leaving little room for improvement by smart repositioning.

\begin{figure}[htbp]
\centering
\includegraphics[width=0.90\textwidth]{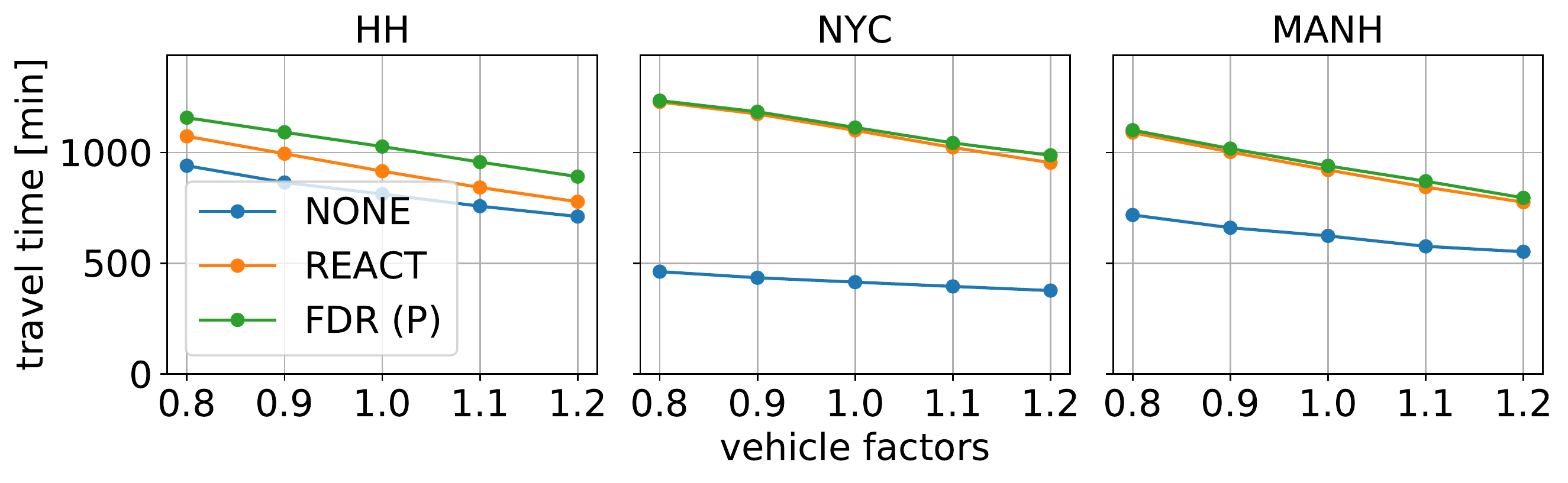}\caption{Average vehicle travel times.}
\label{fig:tt}
\end{figure}

Vehicle travel times and therefore operational costs are increased when using repositioning as seen in Figure \ref{fig:tt}. 
On average, travel times with FDR are increased by $11.7\%$ (HH), $1.7\%$ (NYC) and $2.1\%$ (MANH) compared to REACT. For the HH dataset this increase is larger than might be expected based on the improvement in served requests. One reason for this is that we now also serve those requests that are inefficient to serve, e.g. in remote areas. On the other two datasets the increase is roughly in line with the improvement in served requests.

\subsubsection{Customer waiting times}
\begin{figure}[htbp]
    \centering
    \includegraphics[width=0.90\textwidth]{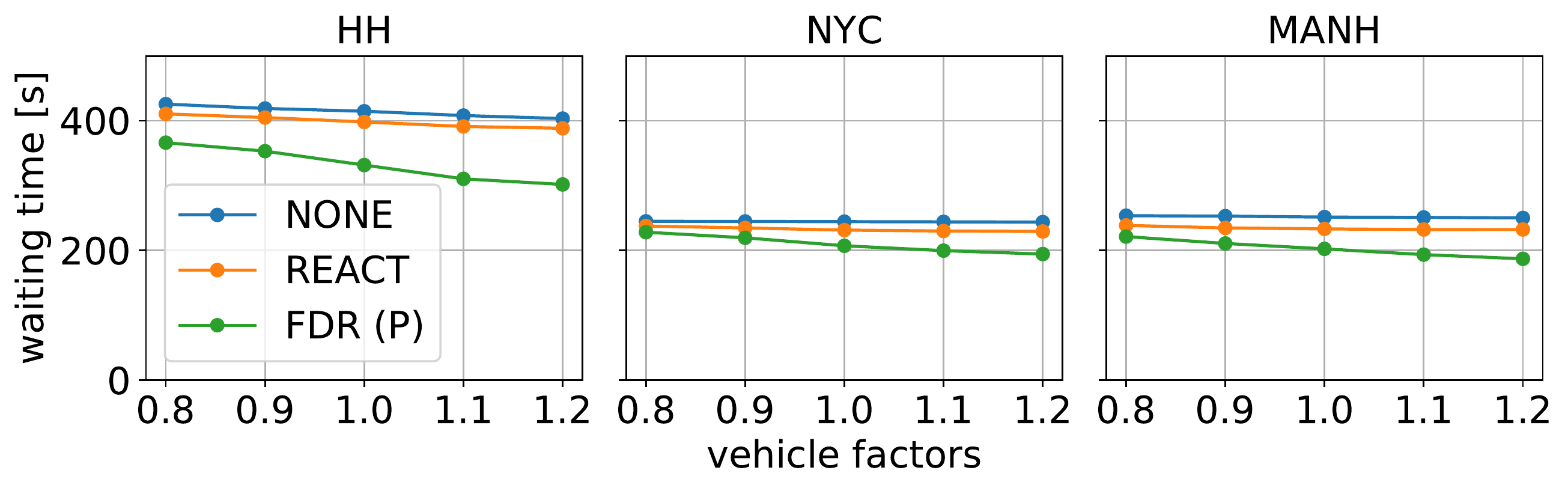}\caption{Average customer waiting times.}
    \label{fig:waiting}
    \end{figure}

Besides reducing rejection rates, repositioning also decreases customer waiting times by moving idle vehicles closer to anticipated customer locations. 
Figure \ref{fig:waiting} compares the average customer waiting time, i.e. the time a customer has to wait before being picked up. 
With FDR the waiting time of a single customer is reduced by an average of $16.6\%$ (66 sec),  $9.9\%$ (23 sec) and $13.2\%$ (31 sec) compared to REACT for HH, NYC, and MANH.
Even in scenarios without a significant improvement regarding rejection rates, FDR manages to reduce waiting times. In practice, this reduction in waiting time will improve customer satisfaction and lead to improved vehicle utilization.

\subsubsection{Vehicle utilization}
\begin{figure}[h]
    \centering
    \includegraphics[width=1.0\textwidth]{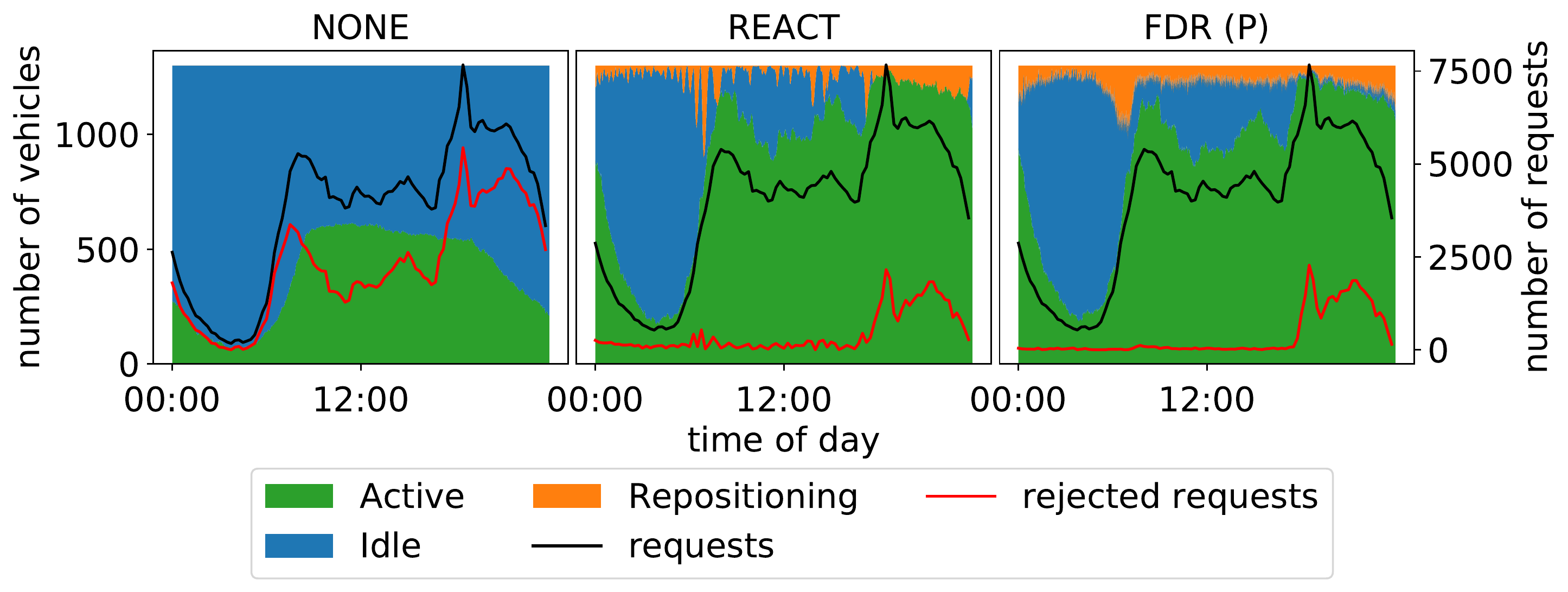}\caption{Vehicle utilization throughout the day for one scenario (NYC, Wednesday, 1.0 vehicle factor) and three repositioning modes. Colored areas illustrate the number of vehicles in a specific state over time. Possible states are idle, active (i.e. serving a tour) and repositioning. Lines indicate the number of total and rejected requests over time}
    \label{fig:utilization}
    \end{figure}

Figure \ref{fig:utilization} shows the vehicle utilization compared between the different repositioning modes for one scenario.
Several aspects of our algorithm FDR may be observed from this chart. During low-demand times (particularly at night between 02:00 and 04:00), most of the fleet is left idle and only minimal repositioning is performed. Before the morning peak, a significant portion of the fleet is repositioned. In comparison, REACT only starts to reposition notable numbers of vehicles after a small spike in rejected requests at around 07:00. Overall, when using FDR, the number of rejected requests is almost zero throughout most of the day. Only during the evening peak after approximately 18:00, when the complete fleet is occupied, requests are rejected.

\subsection{Algorithm performance with a naive demand forecast}
\label{sec:computational_results:naive}
\begin{figure}[htbp]
    \centering
    \includegraphics[width=0.90\textwidth]{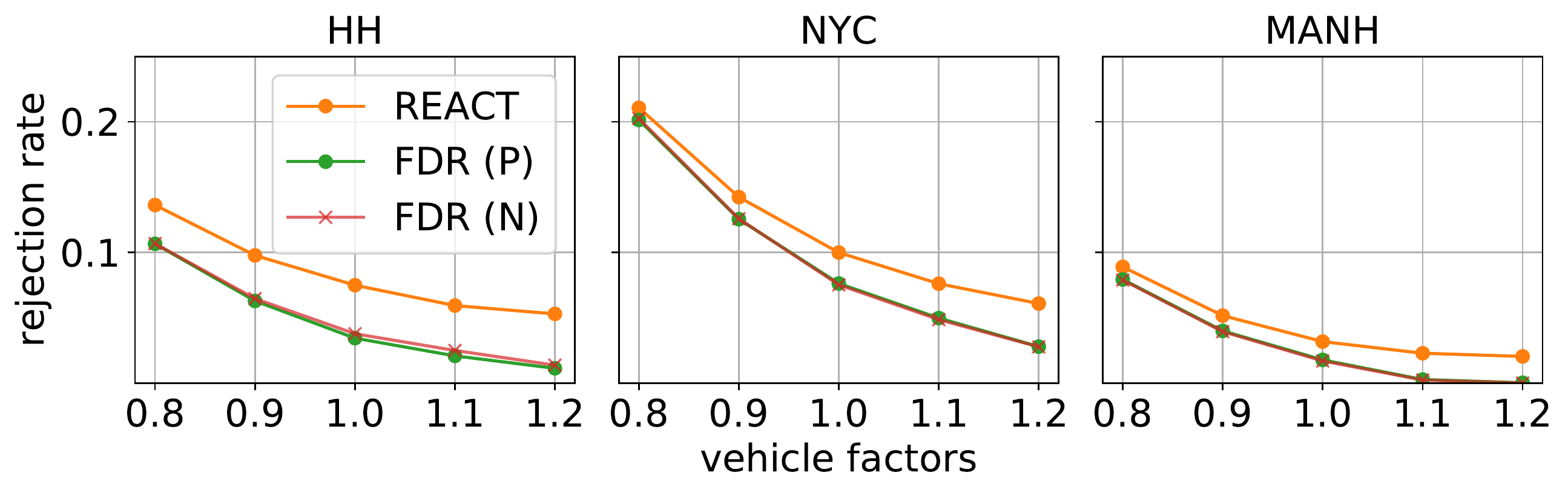}\caption{Average rejection rates.}
    \label{fig:rejection_naive}
    \end{figure}

Figure \ref{fig:rejection_naive} illustrates the average request rejection 
rates with FDR (N) compared to FDR (P) and REACT. 
In comparison with a perfect forecast, the results with a naive forecast are nearly identical with an average increase of 0.06 percentage points.
The picture regarding customer waiting times is similar with an average increase of $0.2 \%$.
These results illustrate that our algorithm is robust to minor forecast errors and may be used successfully with a simple forecasting model. However, it should be noted that for such short-term forecasts the utilized naive model performs remarkably well and would be difficult to outperform substantially even with complex forecasting models.

\section{Conclusion and outlook}
\label{sec:outlook}
In this work, we have presented a forecast-driven algorithm for idle vehicle repositioning. We embedded the algorithm a planning service for dynamic ride-sharing applications and evaluated it through extensive simulations. Our results on three real-world datasets show that our approach is real-time capable even on large-scale scenarios. With a perfect forecast, rejection rates are improved by an average of 2.5 percentage points while customer waiting times are reduced by $13.2\%$. With a naive forecast, results are only slightly worse.

In the future, we aim to study how our algorithm reacts to forecasting errors and in which situations it might lead to undesirable repositioning movements.
Additionally, we intend to improve our model in several ways. The provided coverage of currently traveling vehicles could be modeled in more detail by including spatial-temporal aspects such as current vehicle trajectories. Some model parameters such as the number of trip requests served over the forecast horizon could be determined adaptively and per individual vehicle, increasing the level of detail of the model and removing the need for preliminary parameter optimization.


\end{document}